\documentclass[12pt]{article}

\usepackage{adjustbox}
\usepackage[margin=1in,footskip=0.25in]{geometry}
\usepackage{relsize}
\usepackage{amsthm}%
\usepackage{mathtools}
\usepackage{amsfonts}%
\usepackage{pdfpages}
\usepackage{amsmath}
\usepackage{color}
\usepackage{flexisym}

\usepackage{cite}
\newtheorem{theorem}{Theorem}
\newtheorem{remark}[theorem]{Remark}
\newtheorem{definition}[theorem]{Definition}

\newtheorem{proposition}[theorem]{Proposition}
\newcommand{\CommaBin}{\mathbin{\raisebox{0.5ex}{,}}}

    \DeclareMathOperator{\sech}{sech}

\begin{document}

\title{Hopfield Neuronal Network of Fractional Order: A note on its numerical integration}

\author{Marius-F. Danca{\footnote{Corresponding author}}\\
Romanian Institute of Science and Technology, \\
400487 Cluj-Napoca, Romania\\
Email: danca@rist.ro\\
}

\maketitle

\begin{abstract}
In this paper, the commensurate fractional-order variant of an Hopfield neuronal network is analyzed. The system is integrated with the ABM method for fractional-order equations. Beside the standard stability analysis of equilibria, the divergence of fractional order is proposed to determine the instability of the equilibria. The bifurcation diagrams versus the fractional order, and versus one parameter, reveal a strange phenomenon suggesting that the bifurcation branches generated by initial conditions outside neighborhoods of unstable equilibria are spurious sets although they look similar with those generated by initial conditions close to the equilibria. These spurious sets look ``delayed'' in the considered bifurcation scenario. Once the integration step-size is reduced, the spurious branches maintain their shapes but tend to the branches obtained from initial condition within neighborhoods of equilibria. While the spurious branches move once the integration step size reduces, the branches generated by the initial conditions near the equilibria maintain their positions in the considered bifurcation space. This phenomenon does not depend on the integration-time interval, and repeats in the parameter bifurcation space.
\end{abstract}

\textbf{keyword }Hopfield neural network;
Hidden chaotic attractor; Self-excited attractor; Numerical periodic trajectory; ABM method for fractional-order equations

\vspace{3mm}

\section{Introduction}

Simplified artificial neural networks like the integer-order model considered in \cite{dancus1} and analyzed as a fractional-order model in this paper are inspired by biological neural networks and consists of interconnected groups of neurons.
Chaotic behaviors in these systems are still not adequately addressed nor fully understood (see e.g. \cite{skarda}). The chaos presence in neuronal networks has been investigated in the last thirty years \cite{gid,nara2,cao,aih,free,guck}.

Hopfield Neural Networks (HNNs), particular cases of neural networks, has been inspired from spin systems \cite{spin}. Although if chaos is not easy to identify in these systems, it has been found in many HNNs (see e.g. \cite{csf1,opp,yang,csf3,cpb6}).

In the last thirty years, numerous results on fractional derivative operators suggested that the underlying time memory effect is a useful tool to describe the memory property of an information process. The applications of the Fractional Order (FO) calculus is as old as the Integer Order (IO) one, and was first widely present in mathematics. Because of the nonlocal characteristic of ``infinite memory'' effect, FO systems proved to describe more accurately the behavior of real dynamical systems, compared to the IO models. Definitions of FO derivative have been formulated by Liouville, Grunwald, Letnikov, and Riemann, in the late 19th century, while the first definition of a fractional difference operator was proposed in 1974 \cite{frac20}. For basic aspects of the theory of FO systems see the monographs of Podlubny  \cite{frac10}, Kilbas et al. \cite{kil}.

In this paper, the FO derivative in the sense of Caputo \cite{frac4} is considered especially because it allows the choice of initial conditions as that for the IO systems

\[
D_*^qx=\frac{1}{\Gamma \left \lceil{q}\right \rceil-q }\int_0^t(t-\tau)^{\left \lceil{q}\right \rceil -q-1}D^{\left \lceil{q}\right \rceil }x(\tau)d\tau,
\]
where $D_*^q$ denotes the Caputo differential operator of order $q$ with starting point $0$, $D^{\left \lceil{q}\right \rceil}$ is the standard differential operator of integer order $\left \lceil{q}\right \rceil\in \mathbb{N}$. For $q\in(0,1)$, as considered in this paper, $D^{\left \lceil{q}\right \rceil}x=x\textprime$.
Caputo's derivative models phenomena from the interactions within the past time history (equations having ``memory'') and also for problems with nonlocal properties.

Because the Caputo derivative of a function requires the calculation of its derivative, it is defined only for differentiable functions. Therefore, hereafter it is assumed that all functions are at least differentiable.

The FO HNN of commensurate order considered in this paper belongs to the class of systems modeled by the following autonomous Initial Value Problem

\begin{equation}\label{eqx1}
D_*^qx(t)=f(x(t)),~~~ x(0)=x_0,
\end{equation}
where $f:\mathbb{R}^n\rightarrow \mathbb{R}^n$ is a vector-valued function.

The numerical method used in this paper to integrate the problem \eqref{eqx1} is one of the most performing and utilized numerical methods to integrate systems of FO, the predictor-corrector Adams-Bashforth-Moulton (ABM) method for FDEs \cite{frac3}, which is of a Predict, Evaluate, Correct, Evaluate (PECE) type.

From the computational point of view, there exists the following classification of attractors

\begin{definition}\cite{unu,doi,patru,trei,ixus2}\label{def}
An attractor is called \emph{self-excited} if its basin of attraction intersects
with any open neighborhood of an equilibrium; otherwise, it is called \emph{hidden}.
\end{definition}

Most important ingredients of hidden attractors are: multistability \cite{opt,ixus2}, systems without equilibria \cite{hid1}, systems with stable equilibria \cite{hid1,frac1}, line of equilibria \cite{hid2,hid3}.

While the basin of attraction for a hidden attractor is not connected with
any equilibrium and, therefore, for the numerical localization it is necessary to develop
special analytical-numerical procedures \cite{ixus3,unu,trei}. Self-excited attractors can be visualized numerically by a standard computational procedure, in which
a trajectory starting from a point in a neighborhood of an unstable equilibrium leads to the attractor.

To check that a chaotic attractor of a system which has unstable equilibria is hidden, one verify if the trajectories starting in small neighborhoods of the unstable equilibria are not attracted by the attractor (see e.g. \cite{doi,trei}).

In this paper, we consider the FO variant of the simplified $3$-neuron IO HNN presented in \cite{dancus1}. The FO variant presents some interesting characteristics unveiled by the used Adams-Bashforth-Moulton scheme for FO differential equations.

The paper is structured as follows: Section 2 presents the FO variant of the HNN system and its numerical integration. Section 3 deals with the problem of coexistence of attractors, where the stability of the equilibria and the existence of hidden attractor are studied. In Section 4, a strange dependence on the integration step-size of some attractors is revealed and discussed. The last section of Conclusion closes the paper.

\section{FO Hopfield neural network and its numerical integration}

The FO $3$-neuron HNN considered in this paper is the Caputo FO variant of the simplest example of HNNs of IO analyzed in \cite{dancus1}

\begin{equation}\label{eq}
D_{*}^qx_i=-x_i+\sum_{j=1}^3w_{ij}f(x_j),\quad x_i(0)=x_{0i}\quad i=1,2,3,
\end{equation}

\noindent with the commensurate order $q\in(0,1)$, sigmoid like functions $f(x)=[\tanh(x_1),\tanh(x_2),\tan(x_3)]^t$, approximating the discontinuity at $x_i=0$, with the weight matrix

\[W=
\begin{bmatrix}
    w_{11}& w_{12}& w_{13}\\
    w_{21}& w_{22}& w_{23}\\
    w_{31}& w_{32}& w_{33}
\end{bmatrix}
=
\begin{bmatrix}
1.995&-1.2&0\\
2& 1.71&1.15\\
-4.75&0&1.1
\end{bmatrix}.
\]

Therefore, the system (\ref{eq}) reads
\begin{equation}\label{eq2}
\begin{aligned}
D_*^q{x}_1=&-x_1+1.995~\tanh(x_1)-1.2~\tanh(x_2),\\
D_*^q{x}_2=&-x_2+2~\tanh(x_1)+1.71~\tanh(x_2)+1.15~\tanh(x_3),\\
D_*^q{x}_3=&-x_3-4.75~\tanh(x_1)+1.1~\tanh(x_3).
\end{aligned}
\end{equation}

Chaos, stability analysis, digital implementation or synchronization of HNNs are studied in \cite{Hop_1,Hop_2,Hop_3,Hop_5,Hop_6,Hop_4}.

It is easy to see that the HNN system (\ref{eq2}) is symmetrical with respect to the origin. The equilibria, collinear, are

\[
X_0^*=(0,0,0),\quad X_{1}^*=(0.493,0.366,-3.267),\quad X_{2}^*=(-0.493,-0.366,3.267).
\]

Consider the general system \eqref{eqx1} with the discretization over the numerical time-integration interval $I=[0,T]$, $T>0$, on which the numerical solution is determined, with grid points of an equidistant partition of $I$, $t_i=hi$, $i=0,1,2,...,N$, where $h$ is a fixed step size, $h=T/N$.

First, the method calculates the preliminary approximation $x_{i+1}^P$ for $x(t_{i+1})$ (\emph{predictor phase} determined by the fractional second-order Adams-Bashforth method), as follows

\[
x_{i+1}^P=\sum_{j=0}^{\left \lceil{q}\right \rceil-1}\frac{t_{i+1}^j}{j!}x_0^{(i)}+\frac{1}{\Gamma(q)}\sum_{j=0}^ib_{j,i+1}f(x_j),
\]
where
\[
b_{j,i+1}=\frac{h^q}{q}((i+1)^q-(i-j)^q).
\]
Next, the final approximation of $x(t_{i+1})$, $x_{i+1}$  (\emph{correction phase}, the FO variant of the one-step Adams-Moulton method) is
\[
x_{i+1}=\sum_{j=0}^{\left \lceil{q}\right \rceil-1}\frac{t_{i+1}^j}{j!}x_0^{(i)}+\frac{h^q}{\Gamma{q+2}}\Biggl(\sum_{j=0}^ia_{j,i+1}f(x_j)+f(x_{i+1}^P)\Biggr)
\]
with
\[
a_{j,i+1}=(i-j+2)^{q+1}+(i-j)^{q+1}-2(i-j+1)^{q+1},
\]
for $j=1,2,...,i$.

For $j=0$,
\[
a_{0,i}=i^{q+1}-(i-q)(i+1)^q,
\]
and for $j=i+1$,
\[
a_{i+1,i+1}=1.
\]

Because, $\lceil q\rceil-1=0$, for $q<1$, the first sum in both predictor and corrector phases becomes $x_0$.

Because, at each step, the result is stored for future use in the next integration step, which is why the ABM method has the time history property.

\begin{remark}\label{remix}
The convergence of the method is of order $2$, $\underset{i=1,2,...,N}{\max}|x(t_i)-x_i|=O(h^2)$ \cite{exist}. This convergence is satisfactory from the point of view of the numerical stability. However, the scheme may have a very slow rate of convergence, which deteriorates if $q$ is close to $0$. A solution is to improve the accuracy by introducing supplementary corrector iterations. Improvements of the method can be found in e.g. \cite{exist,frac3} or \cite{frac5};
\end{remark}

\section{Numerical periodic trajectories and attractors coexistence}\label{bifurca}

\begin{figure}
\begin{center}
\includegraphics[scale=0.6]{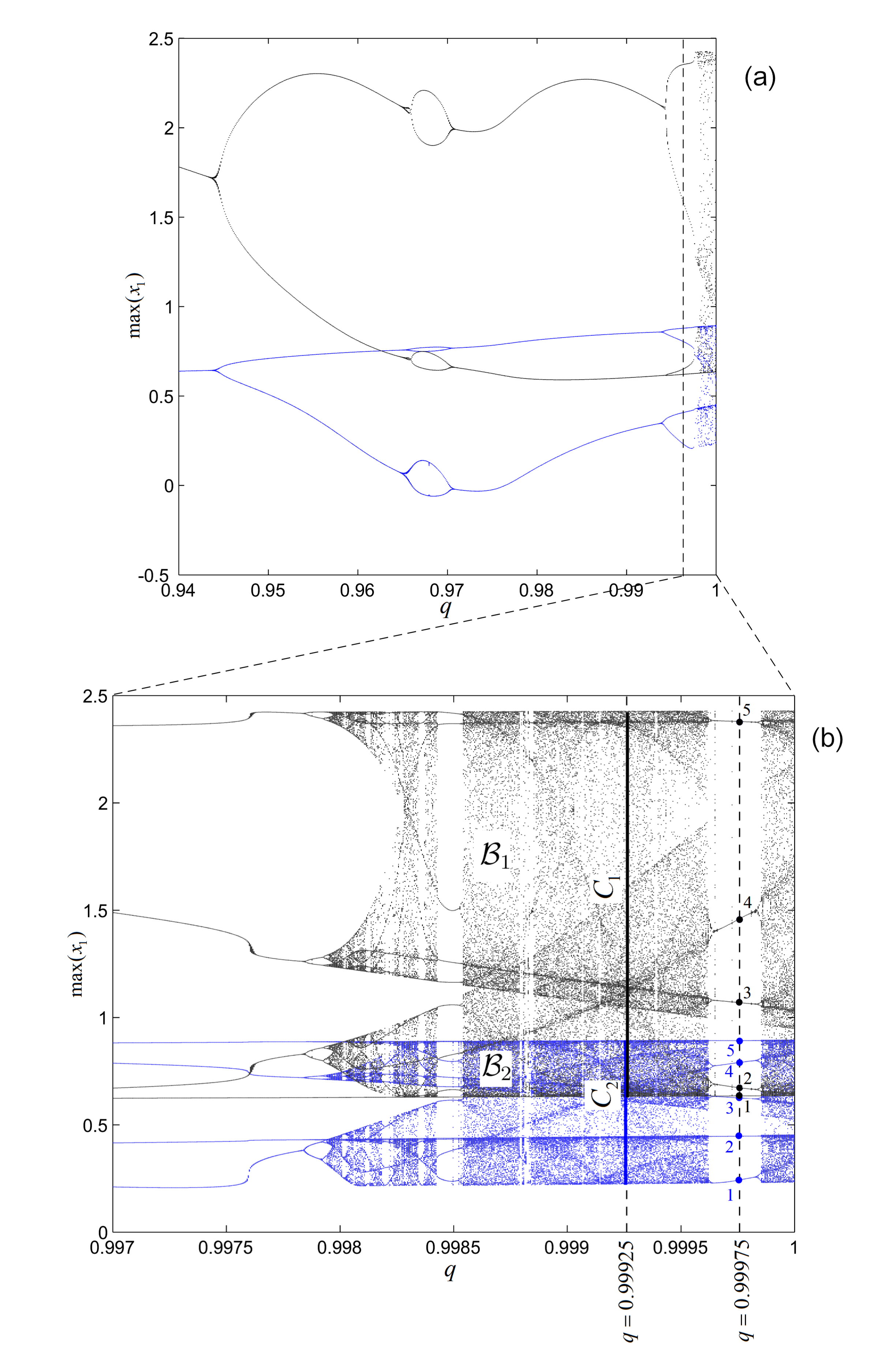}
\caption{Bifurcation diagrams of the system \eqref{eq2} vs the FO $q$; (a) Bifurcation diagram vs the fractional order $q\in[0.94,1]$, for $h=0.01$ and $I=[0,1000]$ and two initial conditions $IC_1=X_1^*$ and $IC_2\approx X_0^*$; (b) Zoom for $q\in[0.997,1]$. The black and blue bifurcative sets, $\mathcal{B}_1$ and $\mathcal{B}_2$, are obtained from $IC_1$ and $IC_2$, respectively. The tick black and blue lines, represent two chaotic behaviors corresponding to $IC_1$ and $IC_2$, respectively, while the filled black and blue circles, the $NPT$s corresponding to $IC_1$ and $IC_2$, respectively.}
\label{fig1}
\end{center}
\end{figure}

One of the best numerical ways to identify the coexistence of attractors is showing the Bifurcation Diagram ($BD$). Due to the system symmetry, the coexistence of symmetric attractors is a feature of the system. Therefore, to simplify the bifurcation pictures of $BD$s, which are symmetric, only the positive local maximum of the state variable $x_1$ is plotted. Fig. \ref{fig1} a) presents the $BD$ vs $q$ for $h=0.01$, $T=1000$ and $q\in[0.94,1]$, while Fig. \ref{fig1} b) for $q\in[0.997,1]$. With these values, $N=T/h=1e5$, which is a relative large value. However, on the other side, for this system the transients for smaller value of $N$ could hide the real attractors (see, for the IO case, \cite{dancus1}).

 Two different branches, $\mathcal{B}_1$ (black color) and $\mathcal{B}_2$ (blue color), are obtained for all considered $BD$s (see Fig. \ref{fig1}). The sets $\mathcal{B}_i$, $i=1,2$, called ``bifurcative'' sets, indicate the coexistence of attractors.
To obtain the $BD$s, the ABM method is applied with the Initial Conditions ($IC$s), $IC_1=X_1^*$, and a close point to $X_0^*$, $IC_2=(1e-3,1e-3,1e-3)$.\footnote{Due to the limited numerical precision, one of the two initial conditions can be chosen $X_1^*$ with three decimals, the second initial condition cannot be chosen as $X_0^*$, but only a close point, because $(0,0,0)$ is an exact solution of the equations which give equilibria.}

 The fact that the rich dynamics occur for higher values of $q$ (close to 1) is typical to continuous-time systems of FO or IO, contrarily to discrete FO systems where phenomena like bifurcation scenario as route to chaos appears for lower values of $q$ (see e.g. \cite{dancus2}).

Cross-sections through $\mathcal{\mathcal{B}}_{1,2}$, on the $BD$, such as the vertical lines at $q=0.99925$ or $q=0.99975$ (Fig. \ref{fig1}), represent linear Poincar\'{e} sections
which
reveal the regular or chaotic behavior of the system for different values of $q$. Therefore, the sets $\mathcal{B}_{1,2}$ can be considered as sets of Poincar\'{e} sections.
For example, if for some $q$ the cross-section contains two dense sets of points (chaotic bands), such as for $q=0.99925$, the system evolves chaotically, while if it contains two finite sets of discrete points (black and blue, for $q=0.99975$), and the system
evolves along a ``stable cycle''.

As known, most dynamics of IO nonlinear systems, refer to stable or unstable periodic solutions (such as limit cycles, chaos, which contains a set of infinitely many unstable periodic orbits mostly with very long periods, quasiperiodic solutions etc). On the other side, in the case of FO systems, the Caputo differential operator cannot transform non-constant periodic functions into periodic ones. Therefore, the following important result was established.

\begin{theorem}\cite{neper1}\label{nu} FO systems modeled by the IVP \eqref{eqx1} cannot have any nonconstant periodic solution.\footnote{Under some circumstances, an FO solution may be asymptotically periodic (see e.g. \cite{neper2}).
}
\end{theorem}

On the other side, many phenomena and real systems are not strictly periodic. Therefore the following definition helps overcome this impediment situation.

\begin{figure}
\begin{center}
\includegraphics[scale=1]{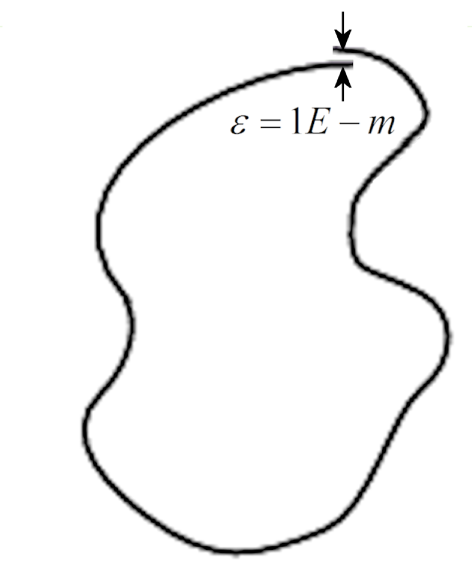}
\caption{Sketch of an $NPT$ in $\mathbb{R}^2$. }
\label{fig2}
\end{center}
\end{figure}

\begin{definition}
In the $n$-dimensional phase space $\mathbb{R}^n$, with $n\geq2$, a \emph{numerically periodic trajectory} ($NPT$) refers to as a closed trajectory in the numerical sense that the closing error $\epsilon$ is within a given bound of $1E -m$, with $m$ being a sufficiently large positive integer (see the sketch in Fig. \ref{fig2}).
\end{definition}

Similarly, since in the case of IO systems the set of unstable periodic trajectories is dense within a chaotic set, from the perspective of Definition \ref{nu}, the \emph{skeleton }of chaos in FO systems can be considered as made from unstable $NPT$s.

\subsection{Stability of equilibria $X_{0,1,2}$}

Consider $q=0.99975$. The cross-section (Fig. \ref{fig1}) indicates the coexistence of two period-$10$ $NPT$s, and the cross-section through both sets $\mathcal{B}_{1,2}$ contains 10 points (see the $2\times5$ filled circles in the BD in Fig. \ref{fig1}): $NPT_1$, corresponding to the $IC_1$, indicated by the black filled circles, and $NPT_2$, corresponding to $IC_2$, indicated by the blue filled circles (see also the phase plot in Figs. \ref{fig3} (b), where the transients are light colored, and the time series are shown in Fig. \ref{fig3} (d)).

The cross-section through $q=0.99925$ (Fig. \ref{fig1}) reveals the coexistence of two chaotic attractors indicated in the cross-section by the black and blue chaotic bands $C_1$ and $C_2$, respectively. The chaotic attractor corresponding to $C_1$, is obtained with $IC_1$, while the attractor corresponding to $C_2$, with $IC_2$ (see the phase plot in Figs. \ref{fig4} (b), and the time series are shown in Fig. \ref{fig4} (d)).

To verify if the system possesses hidden attractors, following Definition \ref{def}, one has to check first the stability of the equilibria.

The Jacobian is\footnote{Because the reciprocal trigonometric functions as the hyperbolic \emph{sech} functions are obscure, Matlab gives an equivalent form for $J$ via \emph{$tanh$} functions.}
\[J=
\begin{bmatrix}
    1.995\sech^2(x_1)-1&-1.2\sech^2(x_2)&0\\
2\sech^2(x_1)&      1.71\sech^2(x_2)-1&1.15\sech^2(x_3)\\
-4.75\sech^2(x_1)&0&1.1\sech^2(x_3)-1
\end{bmatrix},
\]

\noindent and its eigenvalues at $X_0^*$ are $\lambda_1=1.942$ and $\lambda_{2,3}=-0.066 \pm 1.879i$, while the eigenvalues at $X_{1,2}^*$ are $\lambda_1=-0.987$ and $\lambda_{2,3}=0.538 \pm 1.286i$.

Denote the arguments in radians, $\alpha_{i}=\arg(\lambda_i)\in[-\pi,\pi)$ ($angle$ or $atan2$ functions in Matlab), which is the principal branch of the multivalued function Arg$(\lambda)$.

\begin{theorem}Equilibrium $X_0^*$ is unstable for all $q\in(0,1)$.
\begin{proof}
The arguments of eigenvalues are $\alpha_1=0$, and $\alpha_{2,3}=1.6058$.

Consider the quantity \cite{frac1}
\begin{equation}
\iota=q-2\frac{|\alpha_{min}|}{\pi}.
\end{equation}
Because $\alpha_{min}=0$,
\[
\iota=q-2\frac{0}{\pi}=q>0,
\]
which, following the stability theorem for FO systems \cite{frac2,frac30}, implies the instability of the equilibrium $X_0^*$.
\end{proof}
\end{theorem}

By symmetry, for the stability of $X_{1,2}^*$ it is sufficient to only analyze the stability of $X_1$ and, therefore, details are omitted.
\begin{theorem}
Equilibrium $X_{1}^*$ is unstable for $q<0.748$ and stable for $q\geq0.748$.
\begin{proof}
The arguments of eigenvalues are $\alpha_1=\pi$ and $\alpha_{2,3}=1.1745$ and $\alpha_{min}=1.1745$. Therefore
\[
\iota=q-2\frac{1.1745}{\pi}=q-0.7477.
\]
Consequently, for $q<0.748$, $X_1^*$ is stable, while for $q\geq0.748$, $X_1^*$ is unstable.
\end{proof}
\end{theorem}

Analyzing the type and sign of the real component of the three eigenvalues of $X_0^*$, $(+,-,-)$, one can deduce that $X_0^*$ is an attracting focus saddle of index 1 for all $q$ values. As shown by the zoomed rectangle in Fig. \ref{fig3} (b), depending on the attraction basins, trajectories will leave this equilibrium along the unstable manifold of dimension 1 generated by the real positive eigenvalue $\lambda_1$, via spiralling, due to the stable manifold of dimension 2 generated by the negative real part of $\lambda_{2,3}$ (see the zoomed squared around $X_0^*$, where 50 random trajectories are considered).

For $q>0.748$, the case studied in this paper, because of the type and sign of the real components of eigenvalues, $(-,+,+)$, $X_{1,2}^*$ are repelling focus saddles of index 2. Trajectories near these equilibria are attracted along the stable manifold of dimension 1 generated by the real negative eigenvalue $\lambda_1$, and are then rejected via spiraling on the unstable manifold of dimension 2 (due to the positiveness of the real part of $\lambda_{2,3}$ (see the zoomed squared around $X_1^*$ in Fig. \ref{fig3} (f), where 50 random trajectories are considered).

Another possible approach to analyzing the equilibria instability could be by analog with the significance of divergence of IO. Thus, motivated by the results on instability for IO systems presented in \cite{inst1,inst2}, which state that if $Div f(x)>0$ then the point $x$ is unstable, the instability of the equilibria of FO systems could be analyzed via the \emph{fractional divergence}, $Div^{q}$. Therefore, if $Div^{q}f(\bar{x})>0$, the point $\bar{x}$ is unstable.

As known, the Caputo derivative is a linear operator, i.e. for $a,b\in \mathbb{R}$ and $f,g$, some functions such that both $D_*^q f(x)$ and $D_*^q g(x)$ exist
\begin{equation}\label{unu}
D_*^q(af(x)+bg(x))=aD_*^qf(x)+bD_*^qg(x),
\end{equation}
and
\begin{equation}\label{doi}
D_*^q x^n=\frac{\Gamma(n+1)}{\Gamma(n-q+1)}x^{n-q},~ ~\textnormal{for}~~ n\in \mathbb{Z}^+.
\end{equation}
Also
\begin{equation}\label{trei}
D_*^q const=0.
\end{equation}

Let $f=(f_1,f_2,f_3)^t$, $x\mapsto f(x)$, with $x=(x_1,x_2,x_3)\in \mathbb{R}^3$ be a vector valued function, and introduce $Div^q$ in the following form (see also \cite{inst3})

\begin{equation}\label{divgen}
Div^{q} f=\nabla^{q} \cdot f=\bigg(\frac{\partial^q}{\partial x_1^q}\CommaBin \frac{\partial^q}{\partial x_2^q}\CommaBin\frac{\partial^q}{\partial x_3^q}\bigg)\cdot (f_1,f_2,f_3)=\frac{\partial^q f_1}{\partial x_1^q}+\frac{\partial^q f_2}{\partial x_2^q}+\frac{\partial^q f_3}{\partial x_3^q},
\end{equation}
where, by $\frac{\partial^q }{\partial x^q}$, one understands Caputo's derivative with respect to the state variable $x$.

Due to the relative simple form of the Caputo derivative of monomials (see \eqref{doi}), it is obvious that the great majority of nonlinear systems, which are modeled by polynomials, are candidates for this approach.

Because the components of the right-hand function $f$ of the system \eqref{eq2} are $\tanh(x)$, not polynomials, one can approximate them locally with Taylor polynomials.

Consider for simplicity the equilibrium $X_0^*$, for which Taylor approximation of order 5, becomes the following Maclaurin power series, approximated with an sufficiently small error, $O(x^7)$
\[
\tanh(x_i)\approx x_i-\frac{1}{3}x_i^3+\frac{2}{15}x_i^5,\quad i=1,2,3.
\]

Then, using the relations \eqref{unu} and \eqref{doi}, and because Caputo's derivative $\frac{\partial^q}{\partial x_i^q}x_j=0$ for $i\neq j$ (relation \eqref{trei}), one obtains
the derivative of the first component $f_1$ as follows

\begin{equation}
\begin{aligned}
&\frac{\partial^q}{\partial x_1^q}f_1(x)=-\frac{\partial^q}{\partial x_1^q}(x_1)+1.995\frac{\partial^q}{\partial x_1^q}(x_1)-1.995\times\frac{1}{3}\frac{\partial^q}{\partial x_1^q}(x_1^3)+1.995\times\frac{2}{15}\frac{\partial^q}{\partial x_1^q}(x_1^5)-\\
&1.2\underbrace {\frac{\partial^q}{\partial x_1^q}(x_2-\frac{1}{3}x_2^3+\frac{2}{15}x_2^5)}_{\mathlarger{=0}}=-\frac{\Gamma(2)}{\Gamma(2-q)}x_1^{1-q}+1.995\frac{\Gamma(2)}{\Gamma(2-q)}x_1^{1-q}-0.665\frac{\Gamma(4)}{\Gamma(4-q)}x_1^{3-q}+\\
&0.2660\frac{\Gamma(6)}{\Gamma(6-q)}x_1^{5-q}=0.995\frac{\Gamma(2)}{\Gamma(2-q)}x_1^{1-q}-0.665\frac{\Gamma(4)}{\Gamma(4-q)}x_1^{3-q}+0.2660\frac{\Gamma(6)}{\Gamma(6-q)}x_1^{5-q}.
\end{aligned}
\end{equation}
Proceeding similarly, one can obtain the other derivatives, $\frac{\partial^q}{\partial x_i^q}f_i$, $i=2,3$. Finally, via the relation \eqref{divgen}, one obtains $Div^{q}f(x)$ with the graph presented in Fig. \ref{fig5} (blue color), from where one deduces that $Div^{q}f(x)$ at $X_0^*$ is positive for all $q\in(0,1)$ and, therefore, $X_0^*$ is unstable. Similarly, one can deduce the instability for $X_{1,2}^q$, but the computations to approximate $tanh$ within neighborhoods of $X_{1,2}^*$ is too complicated to present here.

Note that, graphically
\begin{equation}\label{div}
\lim_{q\rightarrow 1}Div^q f(x)=Div f(x)=1.805,\quad \text{for}\quad x=(x_1,x_2,x_3)\in V_{X_0^*},
\end{equation}
where $Div f$ is the divergence of the IO system (black color in Fig. \ref{fig5}), which has the form
\[
Div f(x)=-3+1.995\sech^2(x_1)+1.71\sech^2(x_2)+1.1\sech^2(x_3).
\]
Due to the continuity of Caputo's derivative of the right-hand function $f$ of the system \eqref{eq2}, the relation \eqref{div} can also be proved analytically.

\subsection{Hidden and spurious attractors}

\begin{figure}
\begin{center}
\includegraphics[scale=0.55]{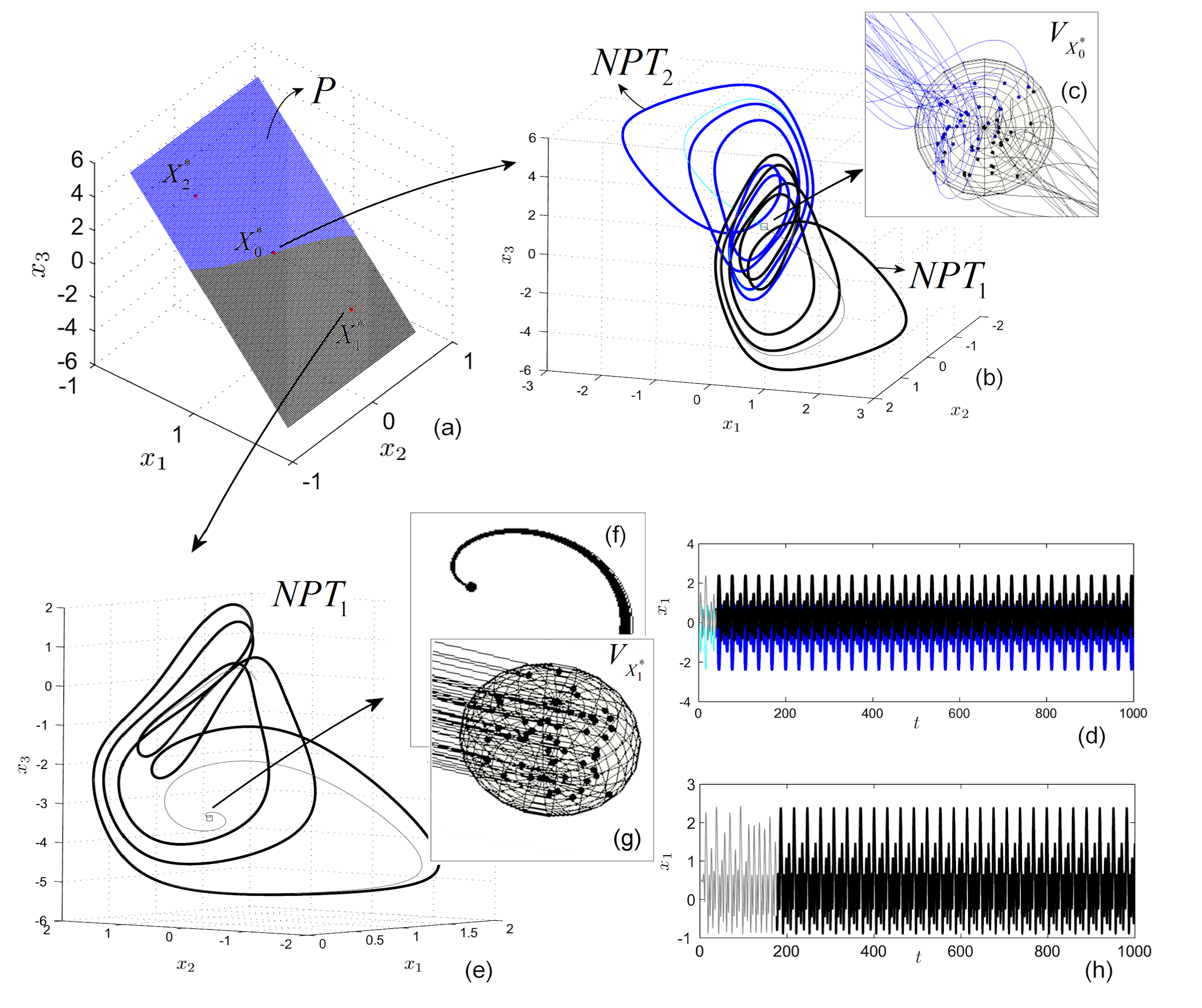}
\caption{$NPT$s for $q=0.99975$. (a) Three-dimensional plane $P$ containing all equilibria $X_{0,1,2}^*$; (b) Phase plots of the two coexisting $NPT$s (transients are light colored); (c) Zoom of the neighborhood of $X_0^*$, within which 50 random $IC$s are considered; (d) Overplot of the time series of the two $NPT$s; (e) Phase plot of the $NPT_1$ generated from a neighborhood of $X_1^*$, $V_{X_1^*}$; (f) The spiralling trajectory starting from $X_1^*$; (g) Zoom of the neighborhood $V_{X_1^*}$, within which 50 random initial conditions are considered; (h) Time series of the first component of $NPT_1$, starting from $X_1^*$.}
\label{fig3}
\end{center}
\end{figure}

\begin{figure}
\begin{center}
\includegraphics[scale=0.45]{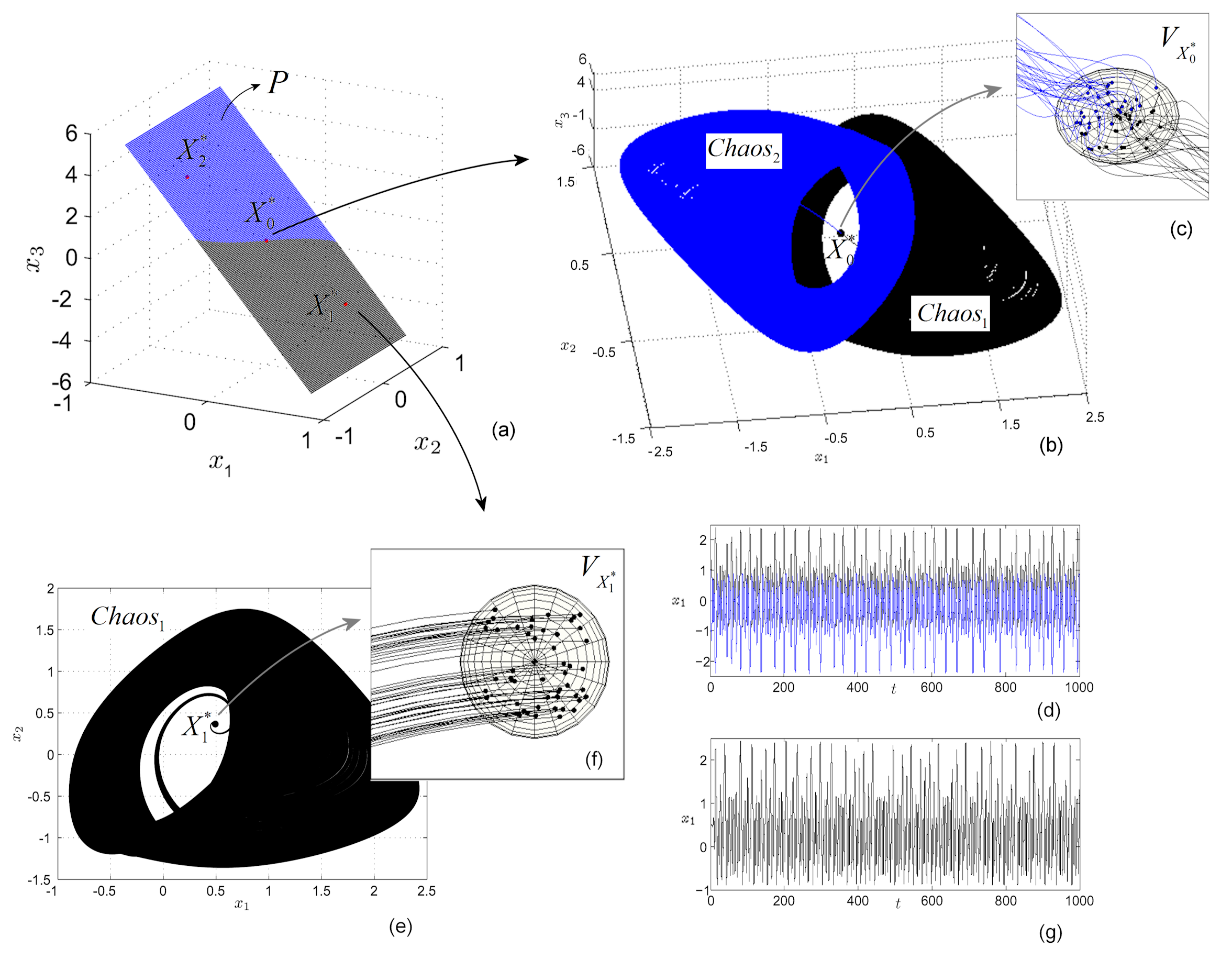}
\caption{Chaotic behavior for $q=0.99925$. (a) Three-dimensional plane $P$ containing all equilibria $X_{0,1,2}^*$; (b) Phase plots of the two coexisting attractors $Chaos_1$ and $Chaos_2$; (c) Zoomed region of the neighborhood of $X_0^*$, within 50 random $IC$s are considered; (d) Overplot of the time series of the two chaotic attractors; (e) Phase plot of the chaotic attractor $Chaos_1$ generated from a neighborhood of $X_1^*$, $V_{X_1^*}$; (f) Zoomed region of the neighborhood $V_{X_1^*}$, within which 50 random initial conditions are considered; (g) Time series of the first component of the chaotic attractor $Chaos_1$, starting from $X_1^*$.
}
\label{fig4}
\end{center}
\end{figure}

\begin{figure}
\begin{center}
\includegraphics[scale=0.65]{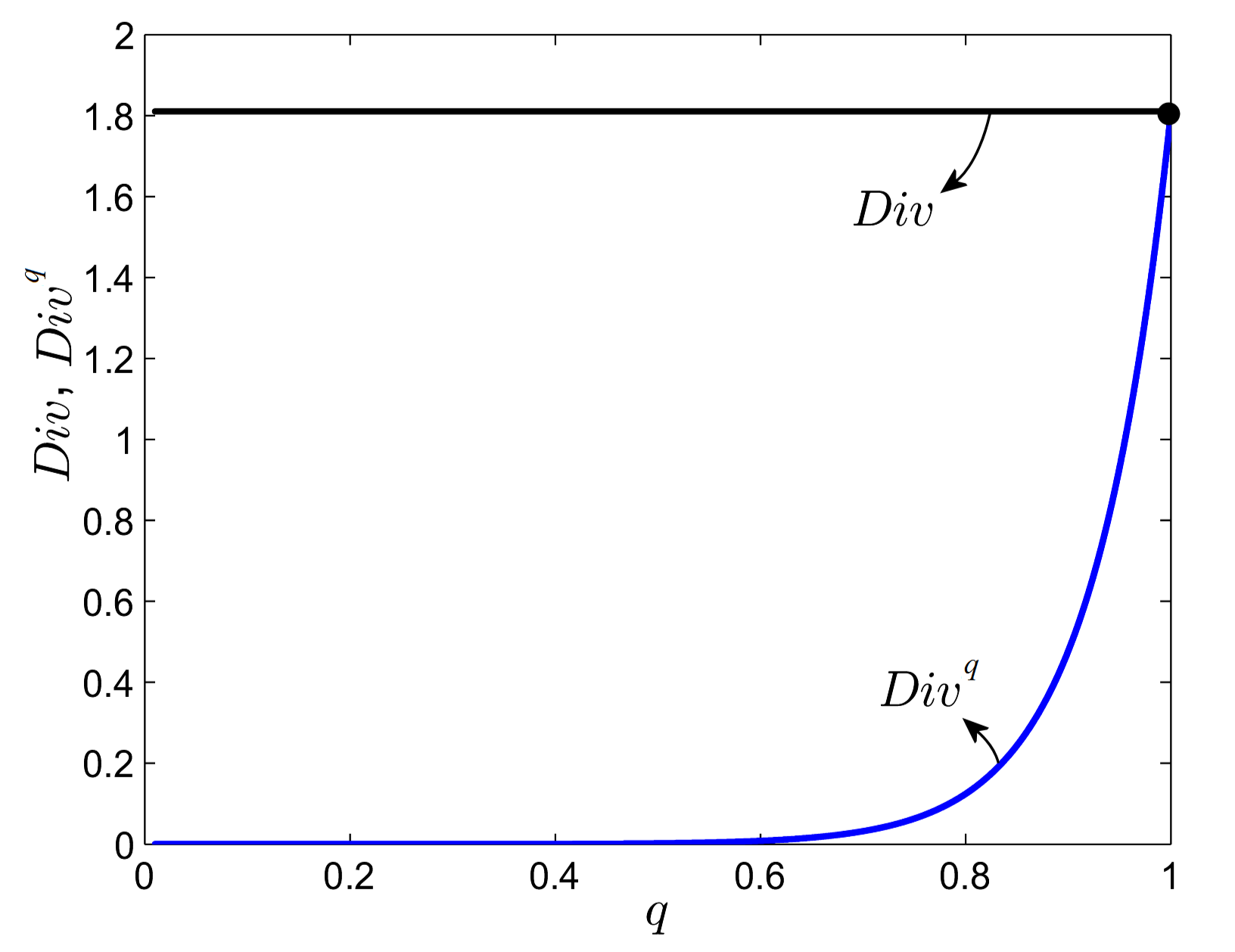}
\caption{Graph of divergence of IO, $Div$ (black), and FO, $Div^q$ (blue), of the system \eqref{eq2}, for $q\in[0,1]$.}
\label{fig5}
\end{center}
\end{figure}

Next, the following result can be numerically verified
\begin{proposition}
The FO HNN system \eqref{eq2} admits no hidden attractors.
\begin{proof}
To verify if an attractor is hidden, usually one takes two-dimensional planar sections through each unstable equilibrium and verify if the trajectories starting within these neighborhoods reach, or not, the attractor. Consider, in this paper, a plane which contains all unstable equilibria $X_{0,1,2}^*$ so that all equilibria can be simultaneously studied. Because the points are collinear, beside the points $X_0^*$, $X_1^*$, consider another non-collinear point, e.g. $(0,1,0)$ on $x_2$ axis. The equation of this plane $P$ is $3.267x_1+0.493x_3=0$, along with a lattice of $300\times 300$ points (Fig.\ref{fig3} (a) and Fig. \ref{fig4} (a)).
\begin{itemize}
\item[i)] Consider the attractors $NPT$s for the case $q=0.99975$ (Fig. \ref{fig1}).
The plane $P$ is scanned, point by point, to find the attraction basins of the two $NPT$s. To each point $(x_1,x_2,x_3)\in P$, considered as $IC$, verify where the underlying trajectories tend. Points $(x_1,x_2,x_3)$, which generate trajectories tending to $NPT_1$, are colored black, while the points generating trajectories tending to $NTP_2$ are colored blue. As can be seen in Fig. \ref{fig3} (a), because all points in neighborhoods of all unstable equilibria tend to one of $NPT_{1,2}$ ($X_0^*$ belongs to the separatrice of attractions basins of $NTP_{1,2}$), these attractors are self-excited (Fig. \ref{fig3} (b)).
 A three-dimensional supplementary numerical verification is presented in the zoomed rectangle in Fig. \ref{fig3} (c) where, within a spherical neighborhood of the equilibrium $X_0^*$, $V_{X_0^*}$, 50 random trajectories are considered. One can see that these trajectories tend to one of the $NPT$s (phase plot in Fig. \ref{fig3} (b) and time series in Fig. \ref{fig3} (d)). A similar three-dimensional approach about $X_1^*$ (Fig. \ref{fig3} (g)) shows that all 50 random trajectories within a small neighborhood of $X_1^*$ tend to $NTP_1$ (see phase plot in Fig. \ref{fig3} (e) and time series in Fig. \ref{fig3} (h), or Fig. \ref{fig3} (a).
\item[ii)] Consider the attractors $Chaos_{1,2}$, for the case of $q=0.99925$, which correspond to the chaotic bands $C_{1,2}$, (Fig. \ref{fig1}). Similarly, like for $NPT$s, scanning the same plane through equilibria $X_{0,1,2}$ (Fig. \ref{fig4} (a)), trajectories tend either to $Chaos_1$ or to $Chaos_2$ (Fig. \ref{fig4} (b)), a fact also underlined by the three-dimensional approach presented in Fig. \ref{fig4} (c). Trajectories within the three-dimensional neoghborhood of $X_1^*$ tend to $Chaos_1$ (Figs. \ref{fig4} (e), (f)).

    Therefore the chaotic attractors $Chaos_{1,2}$ are self-excited.
\end{itemize}
Similarly, it can be shown that the attraction basins of other chaotic attractors and $NPT$s are connected with the unstable equilibria $X_{0,1,2}$ and, therefore, the attractors of this system are self-excited.
\end{proof}
\end{proposition}


The Dahlquist equivalence Theorem provides a tool to check whether or not a numerical scheme is convergent.
\begin{theorem}(Dahlquist Theorem). A multi-step method is convergent if and only if it is consistent and stable.
\end{theorem}
It is said that a multistep numerical method is \emph{consistent }if the underlying discretized equations, for $h\rightarrow 0$, approach the original differential equations. In other words, the consistency defines the relation between the exact solutions and the discrete equations.

A method is \emph{stable} if, in the limit $h\rightarrow 0$, the method has no solutions that grow unbounded as $N=T/h\rightarrow \infty$.

\begin{figure}
\begin{center}
\includegraphics[scale=0.8]{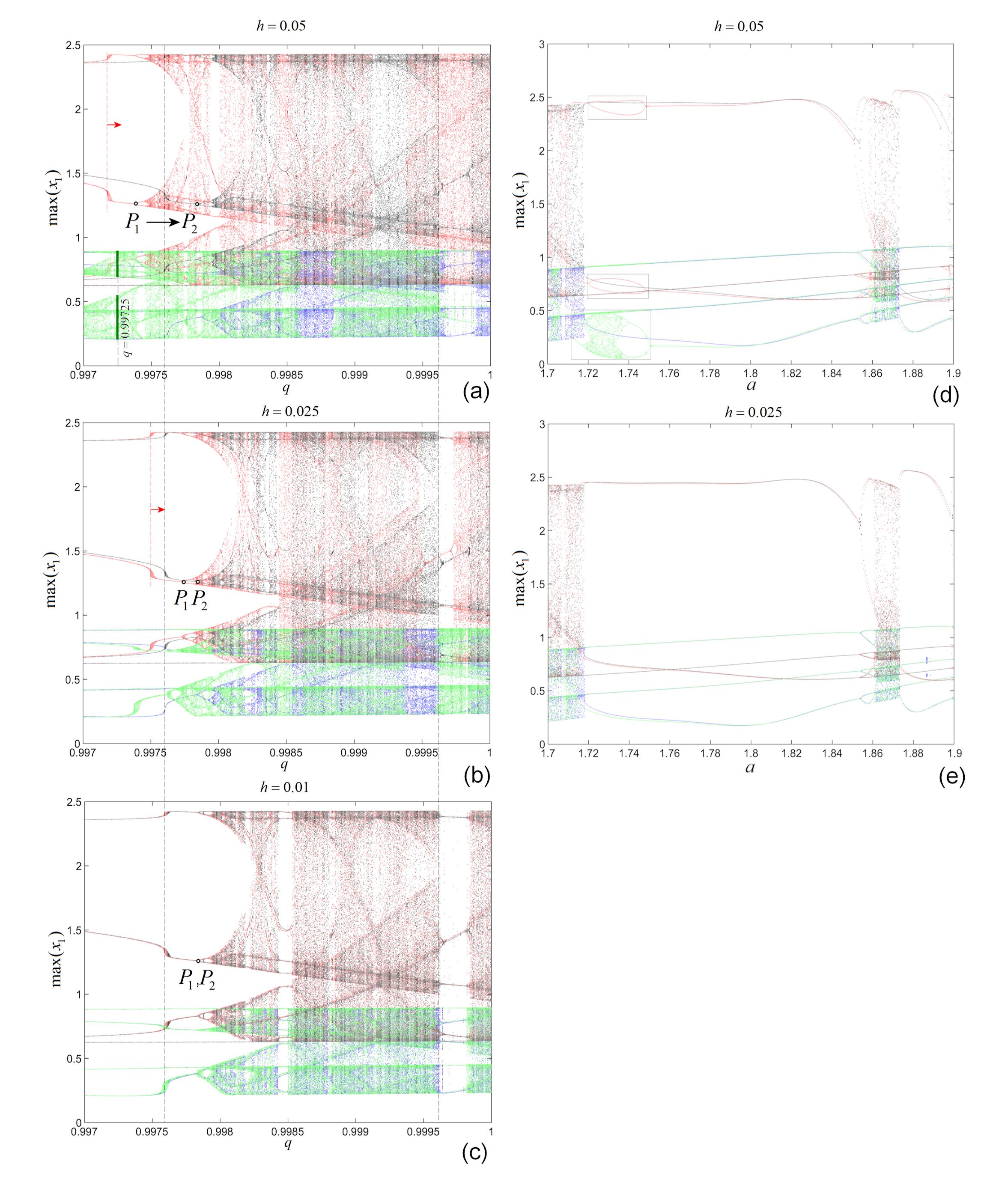}
\caption{Bifurcation diagrams vs $q$ (left column) and vs parameter $w_{11}$ for $q=0.99925$ (right column) for 4 initial conditions; (a)-(c) Bifurcation diagrams vs $q$ for step-size $h=0.05$, $h=0.025$ and $h=0.01$, respectively. Points $P_{1,2}$ indicate the closeness of the spurious bifurcative branches (red and green) to the fixed bifurcative sets (black and blue), once $h$ decreases; (d), (e) Similar phenomenon in the case of bifurcation vs the parameter $w_{11}$.}
\label{fig6}
\end{center}
\end{figure}

As known, the ABM method for FO is convergent (see Remark \ref{remix}). However, in this paper some intriguing phenomenon related to its convergency has been remarked. Fig. \ref{fig6} (left column) presents three $BD$s vs $q$, for three different values of $h$, each of which is obtained with $T=1000$, and four $IC$s: $IC_{1,2}$, which generate the bifurcative sets $\mathcal{B}_{1,2}$ (Fig. \ref{fig1}) and other two $IC$s considered outside the neighborhoods of equilibria.

Consider first the $BD$ for $h=0.05$ (Fig. \ref{fig6} (a)). As can be seen, a false conclusion on the coexistence of attractors could be easily formulated:

\emph{There are four different bifurcative sets}: \emph{the two basic sets, the black and blue, $\mathcal{B}_1$ and $\mathcal{B}_2$ and two new sets, colored red and green}.

However, the bifurcative sets (red and green) are only apparently ``new'', which are spurious sets. Thus, once $h$ reduces, one can see that they tend to the bifurcative sets $\mathcal{B}_1$ and $\mathcal{B}_2$ respectively. Therefore by choosing for example $h=0.025$, the illustrative points $P_1$ on the ``red'' set and $P_2$ on $\mathcal{B}_1$, become closer (compare Fig. \ref{fig6} (a) and Fig. \ref{fig6} (b)), revealing the fact that the new bifurcative sets (red and green) tend to move to right in the bifurcation space with the reduction of $h$. For an even lower value of $h$, $h=0.01$ (Fig. \ref{fig6} (c)), $P_1$ becomes closer to $P_2$, which seems even to coincide with $P_2$ (regardless of inherent computationally errors).

Therefore, with the reduction of $h$, the spurious ``red'' and ``green'' sets tend to the sets $\mathcal{B}_{1,2}$, respectively.

Consequences of this effect is the wrong conclusion that for e.g. $q=0.99725$ due to the spurious chaotic bands (see thick dark green in Fig. \ref{fig6} a)), the system evolves chaotically, or the false conclusion that for this value, chaos and $NPT$s coexist.

\begin{remark}
\begin{itemize}
\item[i)]The reason for this phenomenon is due to the numerical solutions of \eqref{eq2};
\item[ii)] The phenomenon does not depend on the length of the integration time-interval $I$;
\item[iii)] Another interesting fact is that the basic bifurcative sets $\mathcal{B}_{1,2}$, generated from $IC$s belonging to small neighborhoods of unstable equilibria, do not visibly suffer from the numerical effect, but only the bifurcative sets obtained for $IC$s are outside small neighborhoods of equilibria;
\item[iiiv)] The phenomenon repeats for spurious bifurcation sets for whatever supplementary $IC$s are outside neighborhoods from  equilibria;
\item[v)]A similar, even more strange, phenomenon appears on the $BD$ vs parameter $w_{11}\in[1.7,1.9]$ and $q=0.99925$ (the second column in Fig. \ref{fig6}). Considering the same $IC$s and $I$, if for $h=0.05$ the spurious bifurcative sets (``red'' and ``green'') reveal rich dynamics such as bifurcations and chaos (rectangular images in Fig. \ref{fig6} (d)), once $h$ is reduced ($h=0.025$), these dynamics vanish and, moreover, as for the $BD$s vs $q$, but less prominent, the spurious sets tend to the main sets (Fig. \ref{fig6} (e));
\item[vi)] On one side, reducing $h$ implies the vanish of spurious solutions but, on the other side, without some improvements of the ABM scheme (see Remark \ref{remix}), too severe reductions of $h$ (e.g. under $0.01$) imply an extremely long computing time and even possible increase of error \cite{kaik}.
\end{itemize}
\end{remark}

This phenomenon, called here ``$h$-delayed'' because the apparent delay compared to the main bifurcative sets $\mathcal{B}_{1,2}$, could be related more to the consistency of the ABM method, which is strongly related to $h$, rather to the stability which yields the boundedness of solutions.

Concluding, the following proposition could be numerically verified
\begin{proposition}
The results of the numerical integration of the system \eqref{eq2}, with the ABM method, is affected by the $h$-delayed effect.
\end{proposition}

\section{Conclusion}

In this paper, the FO variant of an HNN has been numerically simulated and analyzed. The system is integrated with the ABM method for FDEs. A FO divergence is introduced and used to prove the instability of equilibria. The considered $BD$s versus the FO and versus a system parameter, reveal an interesting phenomenon, named in this paper $h$-delay.
The bifurcative sets generated by initial conditions situated outside equilibria ``move'' in the considered bifurcation space towards the bifurcative sets generated with initial conditions near equilibria, once $h$ reduces. The bifurcative sets generated from equilibria maintain their positions. The phenomenon does not depend on the length of the integration time-interval.
Note that, recently, this ``delay'' phenomenon has been discovered by the author in discrete systems of FO too (paper submitted).
This study will be verified also with other numerical methods for FDEs.

Since the spurious attractors have no connections with equilibria, they can be confused with hidden attractors and, therefore, this phenomenon requires a deep analysis, with determination of an adequate step-size $h$, to avoid false conclusions.

\textbf{Declaration of Competing Interest:} Author declare that he does has no conflict of interest.

\textbf{Funds:}  No funding to declare

\newpage{\pagestyle{empty}\cleardoublepage}

\end{document}